\documentclass{ifacconf}

\usepackage{graphicx}      
\usepackage{natbib}        
\usepackage{amssymb}
\usepackage{amsmath}
\usepackage{dsfont}
\usepackage{euscript}
\usepackage{color}

\newtheorem{definition}{Definition}
\newtheorem{assumption}{Assumption}
\newcommand{\E}{\mathrm{E}}

\newcommand{\R}{\mathrm{Re}}

\newcommand{\changes}[1]{\normalsize{{\color{black} #1}}}
\newcommand{\ch}[1]{\normalsize{{\color{black} #1}}}

\def\R{\mathbb{R}}

\def\xb{\mathbf{x}}

\def\cb{\mathbf{c}}

\def\ab{\mathbf{a}}

\def\yb{\mathbf{y}}
\def\bmu{\boldsymbol{\mu}}

\def\Mb{\boldsymbol{M}}
\def\Qb{\boldsymbol{Q}}
\def\Rb{\boldsymbol{R}}
\def\Ab{\boldsymbol{A}}

\def\R{\mathbb{R}}

\def\xb{\mathbf{x}}

\def\ab{\boldsymbol{a}}
\def\mb{\boldsymbol{m}}
\def\by{\boldsymbol{y}}
\def\cb{\mathbf{c}}

\def\bx{\boldsymbol{x}}
\def\Proj{\mbox{Proj}}
\def\diag{\mbox{diag}}
\def\r{~}

\begin{document}
\begin{frontmatter}

\title{Payoff-Based Approach to Learning Nash Equilibria in Convex Games\thanksref{footnoteinfo}}

\thanks[footnoteinfo]{This research is partially supported by M. Kamgarpour's European Union ERC Starting Grant, CONENE.}

\author[First]{T. Tatarenko}
\author[Second]{M. Kamgarpour}

\address[First]{Department of Control Theory and Robotics, TU Darmstadt, Germany (e-mail: tatarenk@rmr.tu-darmstadt.de)}
\address[Second]{Automatic Control Laboratory, Swiss Federal Institute of Technology, Zurich, Switzerland (e-mail: mkamgar@control.ee.ethz.ch)}

\begin{abstract}                
We consider multi-agent decision making, where each agent optimizes its cost function subject to constraints. Agents' actions belong to a compact convex Euclidean space and the agents' cost functions are coupled. We propose a distributed payoff-based algorithm to learn Nash equilibria in the game between agents. Each agent uses only information about its current cost value to compute its next action. We prove convergence of the proposed algorithm to a Nash equilibrium in the game leveraging established results on stochastic processes. The performance of the algorithm  is analyzed with a numerical case study.
\end{abstract}

\begin{keyword}
Multi-agent decision making, game theory,  payoff-based algorithm
\end{keyword}

\end{frontmatter}


\section{Introduction}
Decision-making in multi-agent systems arises in applications ranging from electricity market to communication and transportation networks \citep{Vehicle, BasharSG, Scutaricdma}. Game theory provides a powerful framework for formulating optimisation problems corresponding to competing or collaborative multi-agent systems. The various notions of equilibria in games  characterise desirable and stable solutions to multi-agent optimisation problems. The focus of our paper is on distributed computation of Nash equilibria for a class of multi-agent decision making modeled by non-cooperative games.

There is a large body of work on computation of Nash equilibria in multi-agent games. The approaches differ by the particular structure of agents' cost functions as well as information available to each agent. In a \textit{potential game},  a central optimization problem can be formulated whose minimizers correspond to Nash equilibria of the game. One can then use  distributed algorithms for computing the minima of the potential function \citep{MardDes, gossipLacra} to converge to Nash equilibria. Distributed algorithms are also proposed for \textit{aggregative games} \citep{Jensen, paccagnan2016aggregative}. In general, for implementation of these distributed algorithms communication is needed between individual agents or between each agent and a central coordinator.

Alternative to distributed optimization approaches, learning approaches to computing Nash equilibria proceed by sampling agents' actions from a set of probability distributions. These probability distributions are updated based on the information available in the system. Most of the past work has focused on algorithms that require knowledge of the structure of the cost functions. For example, \citep{Leslie, MardenLLL, Tat_ACC14, Tat_ECC16} have dealt with learning procedures requiring the so-called oracle information, where each agent can calculate its current cost given any action from its action set.

There are many practical situations in which agents do not know functional form of the objectives. Rather, each agent can only observe their obtained payoffs and be aware of their local actions. In this case, the information structure is referred to as \textit{payoff-based}. A payoff-based learning in potential games is proposed in \citep{MardRev} with the guarantee of stochastic stability of potential function
minimizers, which coincide with Nash equilibria in potential games. However, to implement this payoff-based algorithm agents need to have some memory. Other algorithms requiring only payoff-based information and memory are proposed in \citep{Goto_PIPIP} and \citep{COVEr}. These learning procedures also guarantee stochastic stability of potential function minimizers. Moreover, by tuning a time-dependent parameter the learning procedures converge to a distribution over  potential function minimizers in total variation.

All aforementioned payoff-based procedures are applicable to games with finite action space.
In our recent work, we proposed a payoff-based approach to learn local optima of potential functions in potential continuous action games without memory \citep{Tat_cdc16}. In that work, we addressed the case in which agents' actions live in $\R$.

Our contributions in this paper are as follows. We develop a payoff-based approach for computing Nash equilibria in a class of games with pseudo-monotone maps. In contrast to past work, we consider action spaces being compact subsets of a multidimensional Euclidean space. Given the constraints on action sets in a non-potential game setting, the previously proposed learning methods are no longer applicable. Thus, we develop a sampling based approach, in which the probability distributions from which agents sample their actions are Gaussian. Motivated by projection based algorithms, the mean of the distribution is updated iteratively by each agent based only on its current payoff and projected on the constraint set. This has an interpretation of a stochastic projected gradient algorithm. Using established results on sequence of random variables, and utilizing monotonicity of the game map,  we prove that through appropriate choices of step size, the actions converge in probability to Nash
equilibria.

This paper is organized as follows. In Section \ref{sec:problem} we formulate the game under consideration. In Section \ref{sec:analysis} we present our payoff-based approach and prove its convergence  to a Nash equilibrium. In Section \ref{sec:simulations} we present a numerical case study, motivated by games arising in electricity markets. In Section \ref{sec:conclusion} we summarize the results.

\textbf{Notations and basic definitions:} The set $\{1,\ldots,N\}$ is denoted by $[N]$. Boldface is used to distinguish between the vectors in a multi-dimensional space and scalars.
Given $N$ vectors $\bx^i\in\R^d$, $i\in[N]$, $[\bx^1,\ldots, \bx^N]:=[\bx^i]_{i=1}^{N}:=[{\bx^1}^{\top}, \ldots, {\bx^N}^{\top}]^{\top}$; $\bx^{-i}:=[{\bx^1}, \ldots, {\bx^{i-1}},{\bx^{i+1}}, \ldots, {\bx^N}]$.
The standard inner product on $\R^d$ is denoted by $(\cdot,\cdot)$: $\R^d \times \R^d \to \R$, with associated norm
 $\|\bx\|:=\sqrt{(\bx, \bx)}$. $I_d$ represents the $d$-dimensional identity matrix and $\mathds{1}_N$ represents the $N$-dimensional vector of unit entries. Given some matrix $A\in\R^{d\times d}$, $A\succeq(\succ)0$, if and only if $\bx^{\top}A\bx\ge(>)0$ for all $\bx\ne 0$.
 $A\otimes B$ denotes the Kronecker product. Given a function $g(\bx, \by):\R^{d_1}\times\R^{d_2}\to\R$, we define the
 mapping $\nabla_{\bx}g(\bx, \by): \R^{d_1}\times\R^{d_2}\to\R^{d_1}$ component wise as $[\nabla_{\bx}g(\bx, \by)]_{i}:=\frac{\partial g(\bx, \by)}{\partial x^i}$.
 We will use the big-$O$ notation. Namely, the function $f(x): \R\to\R$ is $O(g(x))$ as $x\to a$, $f(x)$ = $O(g(x))$ as $x\to a$, if $\lim_{x\to a}\frac{|f(x)|}{|g(x)|}\le K$ for some positive $K$.
 \begin{definition}\label{def:pmm}
 The mapping $\Mb:\R^d\to\R^d$ is called \emph{pseudo-monotone} over $X\subseteq\R^d$, if $(\Mb(\by),\bx-\by)\ge 0$ implies $(\Mb(\bx),\bx-\by)\ge 0$ for any $\bx,\by\in X$.
 \end{definition}

  \begin{definition}\label{def:Lip}
 The mapping $\Mb:\R^d\to\R^d$ is called \emph{Lipschitz} with a constant L on $X\subseteq\R^d$, if $\|\Mb(\bx)-\Mb(\by)\|\le L\|\bx-\by\|$ for any $\bx,\by\in X$.
 \end{definition}

\section{Problem Formulation}\label{sec:problem}
\subsection{Learning Nash Equilibria in Games}
We are focused here on a game $\Gamma (N, \{A_i\}, \{J_i\})$ with $N$ players, the sets of players' actions $A_i\subset \R^d$, $i\in[N]$, and the cost functions $J_i:\Ab\to\R$, where $\Ab = A_1\times\ldots\times A_N$ is the set of joint actions\footnote{All results below are applicable for games with different dimensions $\{d_i\}$ of the action sets $\{A_i\}$.}.  We make the following assumptions regarding the game $\Gamma$.

\begin{assumption}\label{assum:convex}
 The game under consideration is \emph{convex}. Namely, for all $i\in[N]$ the set $A_i$ is convex and compact, the cost function $J_i(\ab^i, \ab^{-i})$ is defined on $\R^{Nd}$, continuously differentiable in $\ab$ and convex in $\ab^i$ for  fixed $\ab^{-i}$.
\end{assumption}

\begin{assumption}\label{assum:CG_grad}
 The mapping $\Mb:\R^{Nd}\to\R^{Nd}$, referred to as  the game mapping, defined by
 \begin{align}\label{eq:gamemapping}
 \Mb(\ab) &= [\nabla_{\ab^i} J_i(\ab^i, \ab^{-i})]_{i=1}^N\cr
 &=[M_{1,1}, \ldots, M_{1,d},\ldots, M_{N,1}, \ldots, M_{N,d}]^{\top}, \cr
 M_{i,k}(\ab)&= \frac{\partial J_i(\ab)}{\partial a^i_k}\qquad \ab\in\Ab, \quad i\in[N], \quad k\in[d]
 \end{align}
 is \emph{Lipschitz on $\R^{Nd}$} with a constant $L$ and \emph{pseudo-monotone on $\Ab$} (see Definition \ref{def:pmm}, Definition \ref{def:Lip}).
 \end{assumption}

\begin{assumption}\label{assum:inftybeh}
 The cost functions $J_i(\ab)$, $i\in[N]$, grow not faster than a polynomial function as $\|\ab\|\to\infty$.
\end{assumption}

A Nash equilibrium in a game $\Gamma (N, \{A_i\}, \{J_i\})$ represents a joint action from which no player has any incentive to unilaterally deviate.

\begin{definition}\label{def:NE}
 A point $\ab^*\in\Ab$ is called a \emph{Nash equilibrium} if for any $i\in[N]$ and $\ab^i\in A_i$
 $$J_i(\ab^{i*},\ab^{-i*})\le J_i(\ab^{i},\ab^{-i*}).$$
 \end{definition}
{In this paper, we focus on learning such a stable state in a game. We are interested in designing a \emph{payoff-based algorithm}, which converges to a Nash equilibrium in any game for which Assumptions\r\ref{assum:convex}-\ref{assum:inftybeh} hold.

\subsection{Convex Games and Variational Inequalities}
In this subsection, we prove existence of Nash equilibria for $\Gamma (N, \{A_i\}, \{J_i\})$, through  connecting Nash equilibria and solutions of Variational Inequalities.
\begin{definition}
Consider a mapping $\boldsymbol T(\cdot)$: $\R^d \to \R^d$ and a set $Y \subseteq \R^d$. A
\emph{solution $SOL(Y,\boldsymbol T)$ to the variational inequality problem} $VI(Y,\boldsymbol T)$ is a set of vectors $\yb^* \in Y$ such that $(\boldsymbol T(\yb^*), \yb-\yb^*) \ge 0$, for any $\yb \in Y$.
\end{definition}

The following theorem is the well-known result on the existence of $SOL(Y,\boldsymbol T)$, see Corollary 2.2.5 in \citep{FaccPang1}.
\begin{thm}\label{th:existVI}
 Given $VI(Y,\boldsymbol T)$, suppose that the set $Y$ is compact, convex and that the mapping $\boldsymbol T$ is continuous. Then, $SOL(Y,\boldsymbol T)$ is nonempty and compact.
\end{thm}

Next, we formulate the result that establishes the connection between Nash equilibria in a game and solution vectors of a certain Variational Inequality, Proposition\r1.4.2 in \citep{FaccPang1}.
\begin{thm}\label{th:VINE}
 Given a game $\Gamma (N, \{A_i\}, \{J_i\})$, suppose that the action sets $\{A_i\}$ are closed and convex, the cost functions $\{J_i\}$ are continuously differentiable in $\ab$ and convex in $\ab^i$ for every fixed $\ab^{-i}$ on the interior of $\Ab$. Then, some vector $\ab^*\in \Ab$ is a Nash equilibrium in $\Gamma$, if and only if $\ab^*\in SOL(\Ab,\boldsymbol M)$, where $\Mb$ is the game mapping in \eqref{eq:gamemapping}.
\end{thm}

Note that games for which Assumptions~\ref{assum:convex} and \ref{assum:CG_grad} hold, satisfy all conditions in Theorem\r\ref{th:VINE}. Thus, any solution of $VI(\Ab,\boldsymbol M)$ is also a Nash equilibrium in such games and vice versa. Moreover, according to Theorem\r\ref{th:existVI}, Assumption~\ref{assum:convex} guarantees non-emptiness of $SOL(\Ab,\boldsymbol M)$. Formally, we get the following result.

\begin{cor}\label{cor:existence}
 Let $\Gamma (N, \{A_i\}, \{J_i\})$ be a game for which Assumptions~\ref{assum:convex} and \ref{assum:CG_grad} hold. Then, there exists at least one Nash equilibrium in $\Gamma$. Moreover, any Nash equilibrium in $\Gamma$ belongs to the set $SOL(\Ab,\Mb)$, where $\Mb$ is the game mapping (see \eqref{eq:gamemapping}).
\end{cor}

\ch{
Note that Assumptions~\ref{assum:convex} and \ref{assum:CG_grad} do not imply uniqueness of the Nash equilibrium in $\Gamma (N, \{A_i\}, \{J_i\})$. To guarantee uniqueness, one needs to consider a more restrictive assumption, for example, strong monotonicity of the game mapping \citep{FaccPang1}.  In our paper we do not restrict our attention to such case, but deal with a broader class of games admitting multiple Nash equilibria.}

\section{Solution Approach}\label{sec:analysis}

\subsection{Payoff-Based Algorithm}
In this subsection we formulate the payoff-based approach for the distributed learning of a Nash equilibrium $\ab^*$ in a game $\Gamma(N, \{A_i\}, \{J_i\})$ satisfying Assumptions~\ref{assum:convex}-\ref{assum:CG_grad}.

Having access to information about the current state $\xb^i(t)=[x^i_1,\ldots,x^i_d]^{\top}\in\R^d$ at iteration $t$ and the current cost value $\hat J_i(t)$ at the joint state $\xb(t)$, $\hat J_i(t) = J_i(\xb^1(t),\ldots,\xb^N(t))$, each  agent ``mixes'' its next state $\xb^i(t+1)$, namely it chooses its next state $\xb^i(t+1)$ randomly according to the multidimensional normal distribution $\EuScript N(\bmu^i(t)=[\mu^i_1(t),\ldots,\mu^i_{d}(t)]^{\top},\sigma(t))$ with the density:
\begin{align*}
 p_i&(x^i_1,\ldots,x^i_{d};\bmu^i(t+1),\sigma(t+1))\cr
 & = \frac{1}{(\sqrt{2\pi}\sigma(t+1))^{d}}\exp\left\{-\sum_{k=1}^{d}\frac{(x^i_k-\mu^i_k(t+1))^2}{2\sigma^2(t+1)}\right\}.
\end{align*}
Our choice of Gaussian distribution is based on the idea of CALA (continuous action-set learning automaton), presented in the literature on learning automata \citep{Thatha}.
The mean parameter $\bmu^i(t)$ for the state's distribution is updated as follows:
 \begin{align*}
 &\bmu^i (t+1)=\cr
 & \Proj_{A_i}\left[\bmu^i(t) -\gamma(t+1)\sigma^2(t+1){\hat J_i(t)}\frac{{\xb^i(t)} -\bmu^i(t)}{\sigma^2(t)} \right].
 \end{align*}
In the above, $\Proj_{C}[\cdot]$ denotes the projection operator on set $C$. The initial finite value of $\bmu(0)$ can be defined arbitrarily. We emphasize the difference between states and actions. In particular, states are intermediary values $\xb(t)=[\xb^1(t),\ldots, \xb^N(t)]$ updated during the payoff-based algorithm under consideration. They need not belong to the set of joint actions $\Ab$. We will show that upon convergence of the algorithm, the states will also belong to the joint action set.

Our goal is now to analyze convergence of the proposed algorithm. First, we show that this algorithm is analogous to the Robbins-Monro stochastic approximation procedure \citep{Borkar}. Next, we show the convergence of the random vector $\bmu(t)=[\bmu^1(t),\ldots,\bmu^N(t)]$ using existing results on stochastic processes and by properly choosing $\{\sigma(t), \gamma(t)\}_{t=0}^\infty$.

It is straightforward to show that under Assumption\r\ref{assum:convex}
\begin{align}\label{eq:mathexp}
 &\E_{\xb(t)}\{\hat J_i(t)\frac{x^i_k(t) -\mu^i_k(t)}{\sigma^2(t)}\}  \cr
 =&\E\{\hat J_i(t)\frac{x^i_k(t) -\mu^i_k(t)}{\sigma^2(t)}|x^i_k(t)\sim\EuScript N(\mu_k^i(t),\sigma(t))\}\cr
 = & \frac{\partial {\tilde J_i(\bmu^1(t),\ldots,\bmu^N(t), \sigma(t))}}{\partial \mu^i_k}
\end{align}
 for any $i\in[N]$, $k\in[d]$, where
 \begin{align*}
  \tilde{J}_i &(\bmu^1,\ldots,\bmu^N, \sigma)= \int_{\mathbb R^{Nd}}J_i(\bx)p(\bmu, \bx)d\bx, \cr
  &p(\bmu, \bx)=\prod_{j=1}^Np_j(x^j_1,\ldots,x^j_{d};\bmu^j,\sigma).
 \end{align*}
Note that $\tilde{J}_i$ can be interpreted as the $i$th player's cost function in mixed strategies, given that the mixed strategies are multivariate normal distributions $\{\EuScript N(\bmu^i,\sigma)\}_i$.

We can rewrite the algorithm in the following vector form:
 \begin{align}\label{eq:pbav}
&\bmu(t+1) =\Proj_{\Ab}[\bmu(t) -\gamma(t+1)\sigma^2(t+1)\cr
&\times(\Mb(\bmu(t)) +\Qb(\bmu(t),\sigma(t))+\Rb(\bmu(t),\xb(t),\sigma(t)))],
\end{align}
where
 \begin{align*}
\Qb(&\bmu(t),\sigma(t)) =\tilde{\Mb} (\bmu(t)) -\Mb(\bmu(t)),\cr
\Rb(&\xb(t),\bmu(t),\sigma(t)) = \boldsymbol F(\xb(t),\bmu(t),\sigma(t)) - \tilde{M} (\bmu(t)), \cr
\boldsymbol F(&\xb(t),\bmu(t), \sigma(t)) \cr
&=[{\hat J}_1(t)\frac{\xb^1(t) -\bmu^1(t)}{\sigma^2(t)}, \ldots,  {\hat J}_N(t)\frac{\xb^N(t) -\bmu^N(t)}{\sigma^2(t)}],
\end{align*}
and
$$\tilde{\Mb} (\cdot)=[\tilde M_{1,1}(\cdot), \ldots, \tilde M_{1,d}(\cdot),\ldots, \tilde M_{N,1}(\cdot), \ldots, \tilde M_{N,d}(\cdot)]^{\top}$$
is the $Nd$-dimensional vector, where $i\in[N]$, $k\in[d]$, and
\[\tilde{M}_{i,k} (\bmu(t))=\frac{\partial {\tilde J_i(\bmu^1(t),\ldots,\bmu^N(t), \sigma(t))}}{\partial \mu^i_k}.\]
The algorithm above falls under the framework of Robbins-Monro stochastic approximations procedures  \citep{Borkar}.
Indeed, the vector $\Mb(\bmu(t))$ corresponds to the gradient term in stochastic approximations procedures, $\Qb(\bmu(t),\sigma(t))$ is a disturbance of the gradient term, whereas $\{\Rb(\xb(t), \bmu(t),\sigma(t))\}_t$, according to \eqref{eq:mathexp}, is a martingale difference.
To prove convergence of the algorithm, we will use the following theorem that is a well-known result of Robbins and Siegmund on non-negative random variables, see, for example, Lemma 10 in \citep{Polyak}.

\begin{thm}\label{th:th_nonnegrv} Let $(\Omega, F, P)$ be a probability space and $F_1\subset F_2\subset\dots$ a sequence of sub-$\sigma$-algebras of $F$.
 Let $z_t, \beta_t, \xi_t,$ and $\zeta_t$ be non-negative $F_t$-measurable random variables such that
 \begin{align*}
  \E(z_{t+1}|F_t)\le z_t(1+\beta_t)+\xi_t-\zeta_t.
 \end{align*}
Then almost surely $\lim_{t\to\infty} z_t$ exists and is finite. Moreover, $\sum_{t=1}^{\infty}\zeta_t<\infty$ almost surely on $\{\sum_{t=1}^{\infty}\beta_t<\infty, \sum_{t=1}^{\infty}\xi_t<\infty\}$.
\end{thm}

Now, we are ready to formulate the main result of this section.
\begin{thm}\label{th:main}
  Let players in a game $\Gamma(N, \{A_i\}, \{J_i\})$ update their states $\{\xb^i(t)\}$ at time $t$ according to the normal distribution
  $\EuScript N(\bmu^i(t),\sigma(t))$, where the mean parameters are updated as in \eqref{eq:pbav}. Let Assumptions~\ref{assum:convex}-\ref{assum:inftybeh} hold and the variance parameter $\sigma(t)$ and the step-size parameter $\gamma(t)$ be chosen such that $\sum_{t=0}^{\infty}\gamma(t)\sigma^2(t) = \infty$, $\sum_{t=0}^{\infty}\gamma(t)\sigma^3(t) < \infty$, and
  $\sum_{t=0}^{\infty}\gamma^2(t) < \infty$. Then, as $t\to\infty$, the mean vector $\bmu(t)$ converges almost surely to a Nash equilibrium $\bmu^*\in \Ab$ of the game $\Gamma$,  \ch{given any initial mean vector $\bmu(0)$}, and the joint state $\xb(t)$ converges in probability to $\ab^*=\bmu^*$.
\end{thm}
The theorem above claims almost sure convergence of the sequence of the mean vectors $\{\bmu(t)\}$ and weak convergence of the sequence of the agents' states $\{\xb(t)\}$  to a Nash equilibrium in the game under consideration.

\begin{rem}
 Recall that we distinguish between the states $\{\xb^i\}_{i\in[N]}$ and actions $\{\ab^i\}_{i\in[N]}$ of players in games. During the run of the algorithm, players choose their states $\{\xb^i\}_{i\in[N]}$ according to the normal distributions $\{\EuScript N(\bmu^i,\sigma)\}_{i\in[N]}$ and have access to the current value $\{\hat J_i(t)\}_{i\in[N]}$ of their cost functions, given the actual joint state: $\hat J_i(t) = J_i(\xb^1(t),\ldots,\xb^N(t))$. However, feasibility of the mean vectors $\{\bmu^i(t)\}_{i\in[N]}$ in the proposed procedure justifies the following choice for the actions: $\ab^i(t)=\bmu^i(t)$ for all $i\in[N]$. Thus, under such setting and according to Theorem\r\ref{th:main}, the players' joint action $\ab(t)=[\ab^1(t),\ldots,\ab^N(t)]$ in long run of the payoff-based algorithm converges to a Nash equilibrium almost surely.
\end{rem}

These convergences take place under an appropriate choice of the parameters $\gamma(t)$ and $\sigma(t)$. Note that, analogously to optimization methods based on the gradient descent iterations, the condition $\sum_{t=0}^{\infty}\gamma(t)\sigma^2(t) = \infty$ guarantees sufficient energy for the time-step parameter $\gamma(t)\sigma^2(t)$ to let the algorithm \eqref{eq:pbav} get to a neighborhood of a desired stationary point, \ch{whereas the condition $\sum_{t=0}^{\infty}\gamma^2(t) < \infty$ does not allow the iteration under the projection operator to be unbounded as time goes to infinity.}
\subsection{Proof of Main Result (Theorem \ref{th:main})}
\changes{Our approach is as follows. Firstly, we estimate the distance between the mean vector $\bmu(t+1)$ in the run of the algorithm and some other $\bmu\in\Ab$ by this distance on the previous step, namely by $\|\bmu(t) - \bmu\|$. After that, we analyse each term in this estimation to demonstrate applicability of Theorem\r\ref{th:th_nonnegrv} to the sequence $\{\bmu(t)\}_t$. Finally, we use the properties of Nash equilibria in games satisfying Assumptions~\ref{assum:convex}-\ref{assum:inftybeh} (see Corollary \ref{cor:existence}) to demonstrate that the almost sure limit of the sequence $\{\bmu(t)\}_t$ is a Nash equilibrium in the game under consideration.}

 Let $\beta(t)=\gamma(t)\sigma^2(t)$.
 Then\footnote{We omit further the argument $\sigma(t)$ in terms $\Qb$ and $\Rb$ for the sake of notation simplicity.}
 \begin{align}\label{eq:pbav1}
  \bmu(t+1) &= \Proj_{\Ab}[\bmu(t) -\beta(t+1)\cr
  &\times(\Mb(\bmu(t)) +\Qb(\bmu(t))+\Rb(\bmu(t)))].
 \end{align}

Let $\bmu\in \Ab$ be any point from the joint action set of the game $\Gamma$. Then, taking into account the iterative procedure for the update of $\bmu(t)$ above and the non-expansion property of the projection operator on a convex set, we get
\begin{align}\label{eq:nonexp}
 \|\bmu&(t+1)-\bmu\|^2\cr
 & = \|\Proj_{\Ab}[\bmu(t) -\beta(t+1)\cr
 &\qquad\qquad\times(\Mb(\bmu(t)) +\Qb(\bmu(t))+\Rb(\bmu(t)))]-\bmu\|^2 \cr
 &\le \|\bmu(t) -\beta(t+1)\cr
 &\qquad\qquad\times(\Mb(\bmu(t)) +\Qb(\bmu(t))+\Rb(\bmu(t)))-\bmu\|^2 \cr
 & = \|\bmu(t) - \bmu\|^2 - 2\beta(t+1)(\Mb(\bmu(t)), \bmu(t)-\bmu) \cr
 &\qquad\qquad-2\beta(t+1)(\Qb(\bmu(t))+\Rb(\bmu(t)), \bmu(t)-\bmu) \cr
 &\qquad\qquad + \beta^2(t+1)\|g(\bmu(t))\|^2,
\end{align}
where $g(\bmu(t)) = \Mb(\bmu(t)) +\Qb(\bmu(t))+\Rb(\bmu(t))$.

Let $\EuScript F_T$ be the $\sigma$-algebra generated by the random variables $\{\bmu(k), k\le T\}$.
By taking the conditional expectation with respect to $\EuScript F_T$ of the both sides in the inequality above, we obtain that for any $T>0$, almost surely
\begin{align}\label{eq:condexp}
2\sum_{t=0}^{T}&\beta(t+1)(\Mb(\bmu(t)), \bmu(t)-\bmu)\cr
\le&\|\bmu(0)-\bmu\|^2-\E\{\|\bmu(T+1)-\bmu\|^2|\EuScript F_T\}\cr
+&2\sum_{t=0}^{T}\beta(t+1)\|\Qb(\bmu(t))\|\|\bmu(t)-\bmu)\|\cr
+&\sum_{t=0}^{T}\beta^2(t+1)\E_{\xb(t)}\|g(\bmu(t))\|^2.
\end{align}
In inequality \eqref{eq:condexp} we used the property of the conditional expectation, namely $\E\{\bmu(t_1)|\EuScript F_{t_2}\}=\bmu(t_1)$ almost surely for any $t_1\le t_2$, as well as the fact that $\E\{\Rb(\bmu(t))|\EuScript F_T\}=0$ for all $t\le T$, which is implied by \eqref{eq:mathexp}.

According to Assumption~\ref{assum:inftybeh}, we can show that
\begin{align}\label{eq:diffpotfun}
\tilde{M}_{i,k}(\bmu) &= \frac{1}{\sigma^2}\int_{\R^{Nd}}J_i(\bx)(x^i_k - \mu^i_k)p(\bmu,\bx)d\bx\cr
&=\int_{\R^{Nd}}\frac{\partial J_i(\bx)}{\partial x^i_k}p(\bmu,\bx)d\bx.
 \end{align}
 Thus,
\begin{align}\label{eq:gradmix}
 \tilde{\Mb} (&\bmu(t))=\int_{\mathbb R^{Nd}}{M} (\bx)p(\bmu(t),\bx)d\bx.
\end{align}
Since $\Qb(\bmu(t))=\tilde{\Mb}(\bmu(t)) - \Mb(\bmu(t))$ and due to Assumption~\ref{assum:CG_grad} and equation \eqref{eq:gradmix}, we can write the following:
\begin{align}\label{eq:Qterm1}
 \|\Qb(\bmu(t))\| &\le \int_{\R^{Nd}}\|\Mb(\bmu(t)) - \Mb(\bx)\| p(\bmu(t),\bx) d\bx\cr
 &\le \int_{\R^{Nd}} L \|\bmu(t) - \bx\| p(\bmu(t),\bx) d\bx \cr
 &\le \int_{\R^{Nd}} L \left(\sum_{i=1}^{N}\sum_{k=1}^{d}|\mu^i_k(t) - x^i_k|\right) p(\bmu(t),\bx) d\bx\cr
 &= O(\sigma(t)),
\end{align}
where $L$ is the Lipschitz constant defined in Assumption\r\ref{assum:CG_grad}. The last equality in \eqref{eq:Qterm1} is due to the fact that the first central absolute moment of a random variable with a normal distribution $\EuScript N(\mu,\sigma)$ is $O(\sigma)$.

Obviously, $\|\bmu(t)-\bmu\|$ is bounded for any $t$, since $\bmu(t)\in \Ab$ for any $t$ and $\bmu\in\Ab$.
Now we proceed with estimating the terms $\E_{\xb(t)}\|g(\bmu(t))\|^2$ in \eqref{eq:condexp}:

\begin{align}\label{eq:term2}
\E_{\xb(t)}\|g(\bmu(t))\|^2\le \|\Mb(\bmu(t))\|^2&+\|\Qb(\bmu(t))\|^2\cr
&+2\|\Mb(\bmu(t))\|\|\Qb(\bmu(t))\|\cr
&+\E_{\xb(t)}\|\Rb(\bmu(t))\|^2.
\end{align}

Note that
 \begin{align*}
  \E_{\xb(t)}\|\Rb&(\bmu(t))\|^2\cr
  &\le\sum_{i=1}^N\sum_{k=1}^{d}\int_{\mathbb R^{Nd}}J^2_i(\bx)\frac{(x^i_k - \mu^i_k(t))^2}{\sigma^4(t)}p(\bmu(t),\bx)d\bx.
  \end{align*}
Thus, we can use Assumption~\ref{assum:inftybeh} to conclude that
 \begin{align*}
  \E_{\xb(t)}\|\Rb(\bmu(t))\|^2\le\frac{1}{\sigma^4(t)}f(\bmu(t), \sigma(t)),
 \end{align*}
 where $f(\bmu(t), \sigma(t))$ is a polynomial of $\bmu(t)$ and $\sigma(t)$. Hence, taking into account boundedness of $\bmu(t)$ for all $t$, we conclude that
  \begin{align}\label{eq:term4}
  \beta^2(t+1)\E_{\xb(t)}\|\Rb(\bmu(t))\|^2\le k_1\gamma^2(t)
  \end{align}
for some constant $k_1$.
%
Moreover, according to boundedness of $\bmu(t)$ for all $t$, we can conclude that the first term on the right hand side of \eqref{eq:term2} is bounded.

Bringing \eqref{eq:Qterm1} - \eqref{eq:term4} together and taking into account conditions on the parameters $\gamma(t)$, $\sigma(t)$, we conclude that the right hand side of inequality \eqref{eq:condexp} stays finite almost surely, if $T\to\infty$ and, thus, almost surely
\begin{align}\label{eq:finiteness}
 \sum_{t=0}^{\infty}\beta(t+1)(\Mb(\bmu(t)), \bmu(t)-\bmu)<\infty.
\end{align}

Next, we demonstrate that almost surely
\begin{align}\label{eq:limitpoint}
\varliminf_{t\to\infty} (\boldsymbol M(\bmu(t)), \bmu(t)-\bmu)\le 0.
\end{align}
Indeed, let us assume that, on the contrary, there exists such $\epsilon>0$ and $t_0>0$ that almost surely
$$(\Mb(\bmu(t)), \bmu(t)-\bmu)\ge \epsilon$$
for any $t\ge t_0$. In this case, taking into account that $\sum_{t=0}^{\infty}\beta(t+1)=\infty$, we obtain
\begin{align*}
 \sum_{t=0}^{\infty}\beta(t&+1)(\Mb(\bmu(t)), \bmu(t)-\bmu)\cr
 \ge\sum_{t=0}^{t_0}&\beta(t+1)(\Mb (\bmu(t)), \bmu(t)-\bmu) +\epsilon\sum_{t=t_0}^{\infty}\beta(t+1)=\infty
\end{align*}
almost surely, which contradicts \eqref{eq:finiteness}. Thus, \eqref{eq:limitpoint} holds.
Since $\bmu(t)$ is bounded for any $t$, there exists such a limit point $\bmu^*\in\Ab$ that $\varlimsup_{t\to\infty}\bmu(t)=\bmu^*$ and, according to \eqref{eq:limitpoint},
\begin{align}\label{eq:Nashlimitpoint}
 (\Mb(\bmu^*), \bmu-\bmu^*)\ge 0.
\end{align}
Since we did not specify the choice of $\bmu\in\Ab$, the inequality above holds for any $\bmu\in\Ab$. Thus, according to Corollary\r\ref{cor:existence}, $\bmu^*$ is a Nash equilibrium in the game $\Gamma$.

Next, we notice that, if $\bmu=\bmu^*$ in \eqref{eq:nonexp}, this inequality \eqref{eq:nonexp} together with \eqref{eq:Qterm1} - \eqref{eq:term4} imply that
\begin{align}\label{eq:final}
 \E\{\|\bmu(t+1)-\bmu^*\|^2&|\EuScript F_t\}\le \|\bmu(t)-\bmu^*\|^2\cr
 &-2\beta(t+1)(\Mb(\bmu(t)), \bmu(t)-\bmu^*)\cr
&\qquad\qquad+h(t),
\end{align}
where $h(t) = k_2\beta(t+1)\sigma(t)+k_3\beta^2(t+1)+k_4\beta^2(t+1)\sigma(t)+(\beta(t+1)\sigma(t))^2+k_1\gamma^2(t)$.
According to the properties of $\sigma(t)$ and $\gamma(t)$,
$$\sum_{t=0}^{\infty}h(t)<\infty.$$
Moreover, since $M$ is pseudo-monotone, \eqref{eq:Nashlimitpoint} implies $(\Mb(\bmu(t)), \bmu(t)-\bmu^*)\ge 0$ for any $t$.
Thus, we can apply Theorem\r\ref{th:th_nonnegrv} to conclude that
$$\|\bmu(t)-\bmu^*\| \mbox{ converges almost surely as } t\to\infty.$$
Since $\varlimsup_{t\to\infty}\bmu(t)=\bmu^*$ almost surely,
$$\lim_{t\to\infty}\bmu(t)=\bmu^* \mbox{ almost surely}.$$
\ch{Since $\sum_{t=0}^{\infty}\gamma(t)\sigma^2(t) = \infty$ and $\sum_{t=0}^{\infty}\gamma(t)\sigma^3(t) < \infty$, $\lim_{t\to\infty}\sigma(t)=0$.}
Taking into account that $\xb(t)\sim\EuScript N(\bmu(t),\sigma(t))$, we conclude that $\xb(t)$ converges weakly to a Nash equilibrium $\ab^*=\bmu^*\in\Ab^*$ as time runs. Moreover, according to  Portmanteau Lemma \citep{portlem}, this convergence is also in probability.

\section{Numerical Case Study}\label{sec:simulations}
We illustrate the proposed approach of payoff-based learning to a game arising from electricity market. The agents optimize their power consumption in response to a price signal. The price signal in turn depends on the total consumption of the agents. In contrast to most past work, we consider the  case in which the form of the price function is unknown to agents and the agents do not share information via communication with each other. They can only observe their cost function for every strategy they play.  A similar information structure was considered in \citep{gan2013optimal}.  However, the problem  was decentralized computation of the optimizers of the centralized social welfare objective, rather than a game theoretic formulation.

\subsection{Problem setup}
The problem setup is motivated by the game theoretic formulation of plug-in-electric vehicle (PEV) charging considered in several previous work including \citep{ma2010decentralized,sergio2014tac,couillet2012mean}. The aforementioned problem can be abstracted as a game as follows. There  are $N$ market  participants (users), also referred  to  as players or agents.  Let $\ab^i = [{a}^i_1,\ldots, {a}^i_d]^{\top}\in\mathbb R^d$
denote the decision variable of the player $i$, $i\in[N]$, that is the vector corresponding to its consumption profile over $d$ periods. The constraints for each player $i$ are
\begin{align}\label{eq:constraints}
 0\le &a^i_k\le \bar{a} \quad \mbox{for $k = 1,\ldots,d$},\cr
 \sum_{k=1}^d&a_k^i = \bar{a}^i.
\end{align}
These constraints  indicate that  for  each  player  the
electricity consumption at each time instance is limited and the total electricity consumption over
the considered period of time needs to match a desired amount. The convex and compact set defined by the constraints in \eqref{eq:constraints} is considered the action set $A_i$ for the corresponding player $i$.

The cost function is the price paid for electricity consumption by each agent \citep{paccagnan2016aggregative}
\begin{align}\label{eq:costs}
 J_i(\ab^i,\ab^{-i}) = \ab^{i\top}Q^i\ab^i + 2\left(C^i\frac{1}{N}\sum_{j=1}^N\ab^j+\cb^i\right)^{\top}\ab^i
\end{align}
with $Q^i, C^i\in\mathbb R^{d\times d}$, $\cb^i\in\mathbb R^d$ for all $i\in[N]$. In the above, the first term presents each agent's private value, while the second term corresponds to price of electricity and its functional form may not be known to the agents. 

Consider  the  following  setup.   At  iteration $t$,  each  player  submits  its  proposed  consumption
profile over time horizon of $d$ units, $\xb^i(t) = [{x}^i_1(t), \ldots, {x}^i_d(t)]^{\top}$. How should players update their profiles, using only values of the function $J_i(\xb)$, in order to make the sequence of the joint profiles convergent to a Nash equilibrium? Note that $Q^i$ can in general be known by individual agents, while the second term \eqref{eq:costs} is assumed unknown. Furthermore, $Q^i$ can also equal zero.

For the above game, the game mapping $\Mb(\ab)$ (see \eqref{eq:gamemapping}) is:
\begin{align*}
  \Mb(\ab)= & \hat{M}\ab+\mb,\\
  \hat{M}:= &\diag\left(Q^i+Q^{i\top}+\frac{C^i}{N}\right)+\begin{bmatrix}
           C^{1} \\
           C^{2} \\
           \vdots \\
           C^{N}
         \end{bmatrix}(\mathds{1}_N^{\top}\otimes I_d) \\
  \mb & := [\cb^{1},\; \cb^{2}, \;  \ldots,\;  \cb^{N}].
\end{align*}
We assume the matrices $Q^i$ and $C^i$ to be such that $Q^i+\frac{C^i}{N}\succeq 0$ on $\R^{Nd}$ and $\hat{M}\succeq 0$ on $\Ab$. These conditions imply the fulfillment of Assumption\r\ref{assum:convex}.
Obviously, the mapping $\hat{M}\ab+\mb$ is affine on $\R^{Nd}$ and, hence, Lipschitz on $\R^{Nd}$.
Moreover, the positive semidefinite matrix $\hat{M}$ on $\Ab$ implies that $\hat{M}\ab+\mb$ is pseudo-monotone on $\R^{Nd}$  \citep{AffinePsedoMonot}. Thus, under such setting for the matrices $Q^i$ and $C^i$, $i\in[N]$, Assumptions\r\ref{assum:convex}-\ref{assum:inftybeh} hold in the game under consideration.

\subsection{Simulations}
Let us assume that there are $N$ users, whose strategies are their consumption profiles for $d=4$ periods, the matrices $Q^i$ and $C^i$, $i\in[N]$, in their cost functions \eqref{eq:costs} are the identity matrices of the size $4\times4$, and the vector $\cb^i$, $i\in[N]$, is a $4$-dimensional vector, whose coordinates are some random variables taking values in the interval $(0, 5)$. We assume that the action set $A_i$ for each user $i\in[N]$ is defined by \eqref{eq:constraints}, where $\bar a = 6$ and $\bar{a}_i$ is a random variable taking values in the interval $(0.5,10)$. \ch{The initial mean vector $\bmu(0)$ is a random vector with the uniform distribution on $\Ab=A_1\times\ldots\times A_N$.
Let the agents follow the payoff-based algorithm described by \eqref{eq:pbav}.

Figures~\ref{fig:t10}-\ref{fig:t100} present behavior of the relative error $\frac{\|\bmu(t)-\ab^*\|}{\|\ab^*\|}$ during the algorithm's run with different initial values for $\bmu(0)$, $\gamma(t)=\frac{1}{t^{0.51}}$, $\sigma(t)=\frac{0.1}{t^{0.2}}$, $N=10$ and $N=100$ respectively, where $\ab^*$ is the unique Nash equilibrium of the corresponding game. The uniqueness of the Nash equilibrium is due to the fact that the game mapping is strongly monotone in this example \citep{FaccPang1}.   We can see that already after the first iteration the algorithm gives an approximation for the Nash equilibrium in the game, irrespectively to the initial vector $\bmu(0)$.}
In the case of $N=10$ the relative error is $27\%-31\%$ (see the blue lines on Fig.\r\ref{fig:t10}, Fig.\r\ref{fig:t100}) after the first iteration. In the case of $N=100$, this error drops, however, from $98\%$ and $88\%$ only to $62\%$ and $51\%$ respectively (see the red lines on Fig.\r\ref{fig:t10}, Fig.\r\ref{fig:t100}).
At the next iterations the value of the vector $\bmu(t)$ approaches $\ab^*$ slowly.
After 10 iterations the relative error decreases to $26\%$ and $28\%$, given $N=10$ (see the blue lines on Fig.\r\ref{fig:t10}, Fig.\r\ref{fig:t100}), and to $56\%$ and $45\%$, given $N=100$ (see the red lines on Fig.\r\ref{fig:t10}, Fig.\r\ref{fig:t100}).

From these figures, we observe that convergence of the error to zero is slow. The blue line on Figure\r\ref{fig:t100} demonstrates that the mean vector $\bmu(t)$ is within $85\%$ (the relative error is $15\%$) of the Nash equilibrium in the game with 10 players after 100 iterations, whereas this vector is within $59\%$ (see the red line on Fig.\r\ref{fig:t100}) in the case $N=100$ after 100 iterations. The slow decrease of the relative error after the first iteration can be explained by the choice of the rapidly decreasing parameter $\sigma(t)$ and $\gamma(t)$ as well as by the projection step. These factor prevent a significant change of the projected mean vectors and of the states' values chosen according to the normal distribution with the variance $\sigma(t)$.
Hence, the relation between the settings for the algorithm's parameter, number of agents, and the convergence rate of the procedure should be analyzed in the future.

%
%
\section{Conclusion}\label{sec:conclusion}

\changes{This paper presented a new payoff-based algorithm for learning Nash equilibria in games with pseudo-monotone maps. To investigate the convergence properties of the proposed procedure we used the theory of discrete-time stochastic processes. We proved that in the run of the algorithm the joint actions in the game under consideration converge weakly and in probability to a Nash equilibrium. This payoff-based approach is demonstrated to be applicable to games between users at electrical markets, where the functional form of the electricity price may be unknown to users. Our current and future work focuses on estimation of convergence rate of the algorithm and  improvement of convergence rate by adjustment of the algorithm's parameters.}

\begin{figure}[htb!]
\centering
\includegraphics[width=1\linewidth]{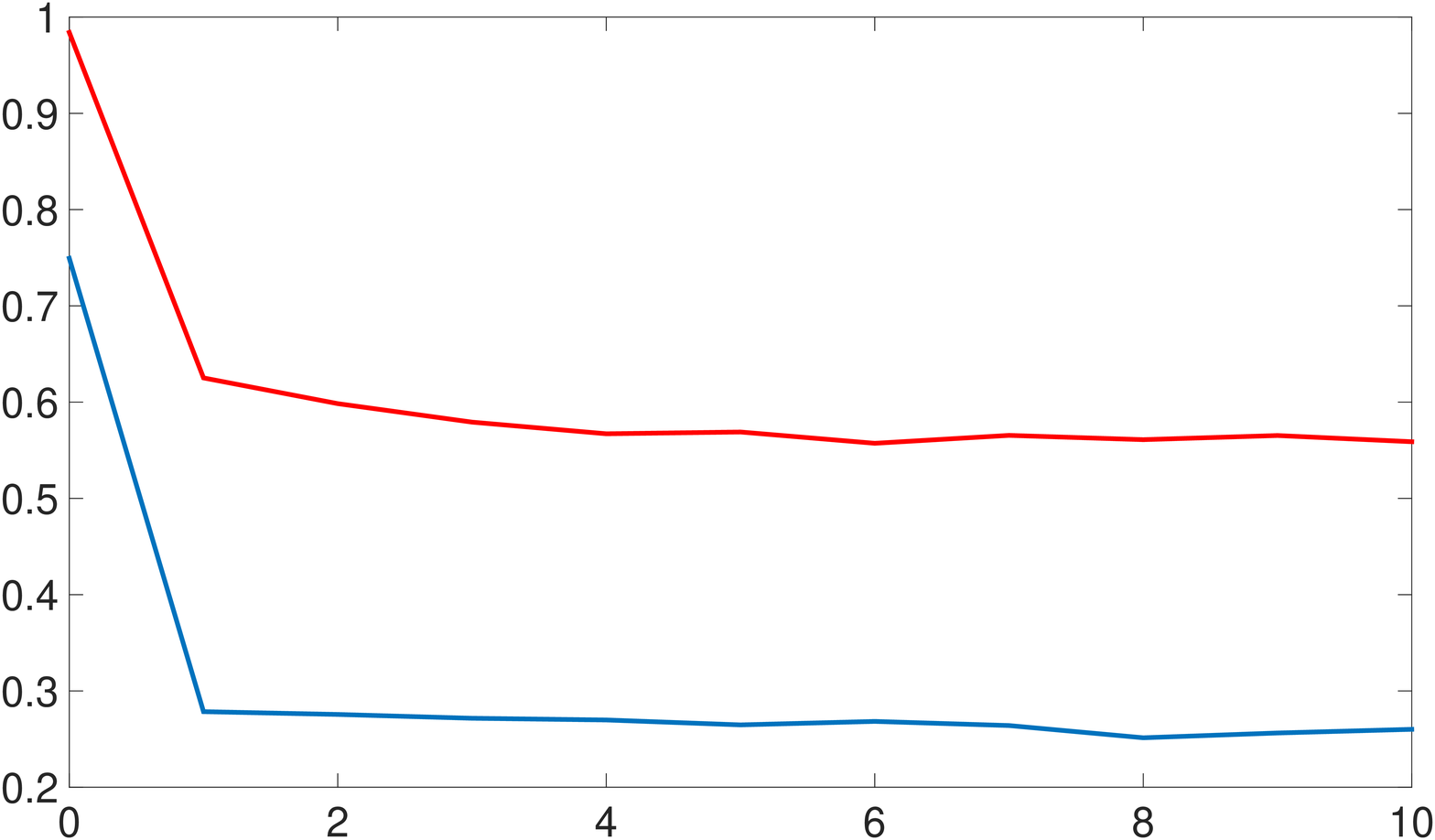}
\caption{Relative error $\frac{\|\bmu(t)-\ab^*\|}{\|\ab^*\|}$ during the payoff-based algorithm, $N=10$ (blue line), $N=100$ (red line), $\gamma(t)=\frac{1}{t^{0.51}}$, $\sigma(t)=\frac{0.1}{t^{0.2}}$.}
\label{fig:t10}
\end{figure}

%

\begin{figure}[htb!]
\centering
\includegraphics[width=1\linewidth]{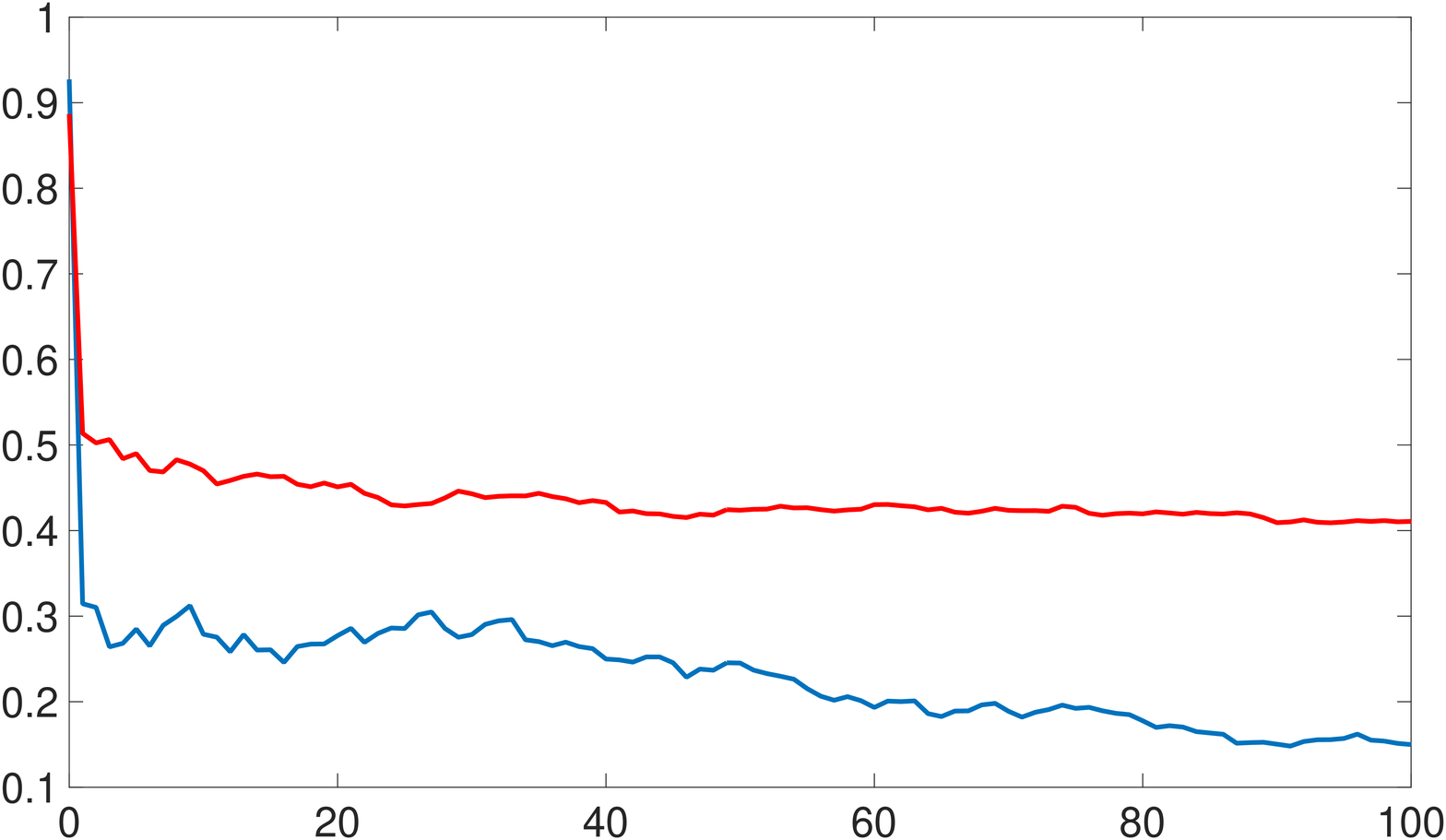}
\caption{Relative error $\frac{\|\bmu(t)-\ab^*\|}{\|\ab^*\|}$ during the payoff-based algorithm, $N=10$ (blue line), $N=100$ (red line), $\gamma(t)=\frac{1}{t^{0.51}}$, $\sigma(t)=\frac{0.1}{t^{0.2}}$.}
\label{fig:t100}
\end{figure}
%
\bibliography{ifacconf}             

\begin{thebibliography}{25}
\providecommand{\natexlab}[1]{#1}
\providecommand{\url}[1]{\texttt{#1}}
\providecommand{\urlprefix}{URL }
\expandafter\ifx\csname urlstyle\endcsname\relax
  \providecommand{\doi}[1]{doi:\discretionary{}{}{}#1}\else
  \providecommand{\doi}{doi:\discretionary{}{}{}\begingroup
  \urlstyle{rm}\Url}\fi

\bibitem[{Arslan et~al.(2007)Arslan, Marden, and Shamma}]{Vehicle}
Arslan, G., Marden, J.R., and Shamma, J.S. (2007).
\newblock Autonomous vehicle-target assignment: a game theoretical formulation.
\newblock \emph{ASME Journal of Dynamic Systems, Measurement and Control}, 129,
  584--596.

\bibitem[{Bharath and Borkar(1999)}]{Borkar}
Bharath, B. and Borkar, V.S. (1999).
\newblock Stochastic approximation algorithms: Overview and recent trends.
\newblock \emph{Sadhana}, 24(4), 425--452.

\bibitem[{Couillet et~al.(2012)Couillet, Perlaza, Tembine, and
  Debbah}]{couillet2012mean}
Couillet, R., Perlaza, S.M., Tembine, H., and Debbah, M. (2012).
\newblock A mean field game analysis of electric vehicles in the smart grid.
\newblock In \emph{{Computer Communications Workshops (INFOCOM WKSHPS)}},
  79--84. IEEE.

\bibitem[{Gan et~al.(2013)Gan, Topcu, and Low}]{gan2013optimal}
Gan, L., Topcu, U., and Low, S.H. (2013).
\newblock Optimal decentralized protocol for electric vehicle charging.
\newblock \emph{IEEE Transactions on Power Systems}, 28(2), 940--951.

\bibitem[{Goto et~al.(2012)Goto, Hatanaka, and Fujita}]{Goto_PIPIP}
Goto, T., Hatanaka, T., and Fujita, M. (2012).
\newblock Payoff-based inhomogeneous partially irrational play for potential
  game theoretic cooperative control: Convergence analysis.
\newblock In \emph{American Control Conference (ACC), 2012}, 2380--2387.

\bibitem[{Gowda(1990)}]{AffinePsedoMonot}
Gowda, M.S. (1990).
\newblock Affine pseudomonotone mappings and the linear complementarity
  problem.
\newblock \emph{SIAM Journal on Matrix Analysis and Applications}, 11(3),
  373--380.

\bibitem[{Grammatico$^*$ et~al.(To Appear Nov. 2016)Grammatico$^*$, Parise$^*$,
  Colombino$^*$, and Lygeros}]{sergio2014tac}
Grammatico$^*$, S., Parise$^*$, F., Colombino$^*$, M., and Lygeros, J. (To
  Appear Nov. 2016).
\newblock Decentralized convergence to {N}ash equilibria in constrained
  deterministic mean field control.
\newblock \emph{\textup{accepted to} IEEE Transactions on Automatic Control}.
\newblock \urlprefix\url{http://arxiv.org/pdf/1410.4421.pdf}.

\bibitem[{Jensen(2010)}]{Jensen}
Jensen, M.K. (2010).
\newblock Aggregative games and best-reply potentials.
\newblock \emph{Economic Theory}, 43(1), 45--66.

\bibitem[{Klenke(2008)}]{portlem}
Klenke, A. (2008).
\newblock \emph{Probability theory: a comprehensive course}.
\newblock Springer, London.

\bibitem[{Li and Marden(2013)}]{MardDes}
Li, N. and Marden, J.R. (2013).
\newblock Designing games for distributed optimization.
\newblock \emph{IEEE Journal of Selected Topics in Signal Processing}, 7(2),
  230--242.
\newblock Special issue on adaptation and learning over complex networks.

\bibitem[{Ma et~al.(2010)Ma, Callaway, and Hiskens}]{ma2010decentralized}
Ma, Z., Callaway, D., and Hiskens, I. (2010).
\newblock Decentralized charging control for large populations of plug-in
  electric vehicles.
\newblock In \emph{49th IEEE conference on decision and control (CDC)},
  206--212. IEEE.

\bibitem[{Marden et~al.(2009)Marden, Arslan, and Shamma}]{MardenLLL}
Marden, J.R., Arslan, G., and Shamma, J.S. (2009).
\newblock Cooperative control and potential games.
\newblock \emph{Trans. Sys. Man Cyber. Part B}, 39(6), 1393--1407.

\bibitem[{Marden and Shamma(2012)}]{MardRev}
Marden, J.R. and Shamma, J.S. (2012).
\newblock Revisiting log-linear learning: Asynchrony, completeness and
  payoff-based implementation.
\newblock \emph{Games and Economic Behavior}, 75(2), 788 -- 808.

\bibitem[{Paccagnan et~al.(2016)Paccagnan, Kamgarpour, and
  Lygeros}]{paccagnan2016aggregative}
Paccagnan, D., Kamgarpour, M., and Lygeros, J. (2016).
\newblock {On Aggregative and Mean Field Games with Applications to Electricity
  Markets}.
\newblock In \emph{European Control Conference}.

\bibitem[{Pang and Facchinei(2003)}]{FaccPang1}
Pang, J.S. and Facchinei, F. (2003).
\newblock \emph{Finite-dimensional variational inequalities and complementarity
  problems : vol. 2}.
\newblock Springer series in operations research. Springer, New York, Berlin,
  Heidelberg.

\bibitem[{Perkins et~al.(2015)Perkins, Mertikopoulos, and Leslie}]{Leslie}
Perkins, S., Mertikopoulos, P., and Leslie, D.S. (2015).
\newblock Mixed-strategy learning with continuous action sets.
\newblock \emph{IEEE Transactions on Automatic Control}, (open access).

\bibitem[{Poljak(1987)}]{Polyak}
Poljak, B.T. (1987).
\newblock \emph{Introduction to optimization}.
\newblock Optimization Software.

\bibitem[{Saad et~al.(2012)Saad, Zhu, Poor, and Basar}]{BasharSG}
Saad, W., Zhu, H., Poor, H.V., and Basar, T. (2012).
\newblock Game-theoretic methods for the smart grid: An overview of microgrid
  systems, demand-side management, and smart grid communications.
\newblock \emph{IEEE Signal Processing Magazine}, 29(5), 86--105.

\bibitem[{Salehisadaghiani and Pavel(2014)}]{gossipLacra}
Salehisadaghiani, F. and Pavel, L. (2014).
\newblock Nash equilibrium seeking by a gossip-based algorithm.
\newblock In \emph{53rd IEEE Conference on Decision and Control}, 1155--1160.

\bibitem[{Scutari et~al.(2006)Scutari, Barbarossa, and Palomar}]{Scutaricdma}
Scutari, G., Barbarossa, S., and Palomar, D.P. (2006).
\newblock Potential games: A framework for vector power control problems with
  coupled constraints.
\newblock In \emph{2006 IEEE International Conference on Acoustics Speech and
  Signal Processing Proceedings}, volume~4, 241--244.

\bibitem[{Tatarenko(2014)}]{Tat_ACC14}
Tatarenko, T. (2014).
\newblock Proving convergence of log-linear learning in potential games.
\newblock In \emph{American Control Conference (ACC), 2014}, 972--977.

\bibitem[{Tatarenko(2016{\natexlab{a}})}]{Tat_cdc16}
Tatarenko, T. (2016{\natexlab{a}}).
\newblock Stochastic payoff-based learning in multi-agent systems modeled by
  means of potential games.
\newblock In \emph{55th IEEE Conference on Decision and Control}. accepted.

\bibitem[{Tatarenko(2016{\natexlab{b}})}]{Tat_ECC16}
Tatarenko, T. (2016{\natexlab{b}}).
\newblock Stochastic stability of potential function maximizers in continuous
  version of independent log-linear learning.
\newblock In \emph{European Control Conference (ECC), 2016}.

\bibitem[{Thathachar and Sastry(2003)}]{Thatha}
Thathachar, A.L. and Sastry, P.S. (2003).
\newblock \emph{Networks of Learning Automata: Techniques for Online Stochastic
  Optimization}.
\newblock Springer US.

\bibitem[{Zhu and Mart\'{\i}nez(2013)}]{COVEr}
Zhu, M. and Mart\'{\i}nez, S. (2013).
\newblock Distributed coverage games for energy-aware mobile sensor networks.
\newblock \emph{SIAM J. Control and Optimization}, 51(1), 1--27.

\end{thebibliography}








\end{document}